\documentstyle[epic,eepic]{article}
\def\hybrid{\topmargin 0pt      \oddsidemargin 0pt
        \headheight 0pt \headsep 0pt
        \textwidth 160true mm       
        \textheight 231true mm         
        \marginparwidth 0.0in
        \parskip 0pt plus 1pt   \jot = 1.5ex}
\catcode`\@=11
\def\marginnote#1{}

\newcount\hour
\newcount\minute
\newtoks\amorpm
\hour=\time\divide\hour by60
\minute=\time{\multiply\hour by60 \global\advance\minute by-\hour}
\edef\standardtime{{\ifnum\hour<12 \global\amorpm={am}%
        \else\global\amorpm={pm}\advance\hour by-12 \fi
        \ifnum\hour=0 \hour=12 \fi
        \number\hour:\ifnum\minute<10 0\fi\number\minute\the\amorpm}}
\edef\militarytime{\number\hour:\ifnum\minute<10 0\fi\number\minute}

\def\draftlabel#1{{\@bsphack\if@filesw {\let\thepage\relax
   \xdef\@gtempa{\write\@auxout{\string
      \newlabel{#1}{{\@currentlabel}{\thepage}}}}}\@gtempa
   \if@nobreak \ifvmode\nobreak\fi\fi\fi\@esphack}
        \gdef\@eqnlabel{#1}}
\def\@eqnlabel{}
\def\@vacuum{}
\def\draftmarginnote#1{\marginpar{\raggedright\scriptsize\tt#1}}

\def\draft{\oddsidemargin -.5truein
        \def\@oddfoot{\sl preliminary draft \hfil
        \rm\thepage\hfil\sl\today\quad\militarytime}
        \let\@evenfoot\@oddfoot \overfullrule 3pt
        \let\label=\draftlabel
        \let\marginnote=\draftmarginnote
   \def\@eqnnum{(\theequation)\rlap{\kern\marginparsep\tt\@eqnlabel}%
\global\let\@eqnlabel\@vacuum}  }
\newcounter{app}
\newcounter{sapp}[app]

\makeatletter
\newdimen\normalarrayskip              
\newdimen\minarrayskip                 
\normalarrayskip\baselineskip
\minarrayskip\jot
\newif\ifold             \oldtrue            
\def\arraymode{\ifold\relax\else\displaystyle\fi} 
\def\eqnumphantom{\phantom{\mbox{\rm
(\theequation)}}}
\def\@arrayskip{\ifold\baselineskip\z@\lineskip\z@
     \else
     \baselineskip\minarrayskip\lineskip2\minarrayskip\fi}
\def\@arrayclassz{\ifcase \@lastchclass \@acolampacol \or
\@ampacol \or \or \or \@addamp \or
   \@acolampacol \or \@firstampfalse \@acol \fi
\edef\@preamble{\@preamble
  \ifcase \@chnum
     \hfil$\relax\arraymode\@sharp$\hfil
     \or $\relax\arraymode\@sharp$\hfil
     \or \hfil$\relax\arraymode\@sharp$\fi}}
\def\@array[#1]#2{\setbox\@arstrutbox=\hbox{\vrule
     height\arraystretch \ht\strutbox
     depth\arraystretch \dp\strutbox
     width\z@}\@mkpream{#2}\edef\@preamble{\halign \noexpand\@halignto
\bgroup \tabskip\z@ \@arstrut \@preamble \tabskip\z@ \cr}%
\let\@startpbox\@@startpbox \let\@endpbox\@@endpbox
  \if #1t\vtop \else \if#1b\vbox \else \vcenter \fi\fi
  \bgroup \let\par\relax
  \let\@sharp##\let\protect\relax
  \@arrayskip\@preamble}
\def\eqnarray{\stepcounter{equation}%
              \let\@currentlabel=\theequation
              \global\@eqnswtrue
              \global\@eqcnt\z@
              \tabskip\@centering
              \let\\=\@eqncr
              $$%
 \halign to \displaywidth\bgroup
    \eqnumphantom\@eqnsel\hskip\@centering
    $\displaystyle \tabskip\z@ {##}$%
    &\global\@eqcnt\@ne \hskip 1.2\arraycolsep
         $\displaystyle\arraymode{##}$\hfil
    &\global\@eqcnt\tw@ \hskip 1.2\arraycolsep
         $\displaystyle\tabskip\z@{##}$\hfil
         \tabskip\@centering
    &{##}\tabskip\z@\cr}
\makeatother
\catcode`\@=11
\ifcase\@ptsize
 \font\tenmsa=msam10
 \font\sevenmsa=msam7
 \font\fivemsa=msam5
 \font\tenmsb=msbm10
 \font\sevenmsb=msbm7
 \font\fivemsb=msbm5
 \font\teneu=eufm10
 \font\seveneu=eufm7
 \font\fiveeu=eufm5
 \font\tenib=cmmib10
 \font\sevenib=cmmib7
 \font\fiveib=cmmib5
\or
 \font\tenmsa=msam10 scaled \magstephalf
 \font\sevenmsa=msam7 scaled \magstephalf
 \font\fivemsa=msam5 scaled \magstephalf
 \font\tenmsb=msbm10 scaled \magstephalf
 \font\sevenmsb=msbm7 scaled \magstephalf
 \font\fivemsb=msbm5  scaled \magstephalf
 \font\teneu=eufm10  scaled \magstephalf
 \font\seveneu=eufm7  scaled \magstephalf
 \font\fiveeu=eufm5   scaled \magstephalf
 \font\tenib=cmmib10  scaled \magstephalf
 \font\sevenib=cmmib7  scaled \magstephalf
 \font\fiveib=cmmib5   scaled \magstephalf
\or
 \font\tenmsa=msam10 scaled \magstep1
 \font\sevenmsa=msam7 scaled \magstep1
 \font\fivemsa=msam5  scaled \magstep1
 \font\tenmsb=msbm10 scaled \magstep1
 \font\sevenmsb=msbm7 scaled \magstep1
 \font\fivemsb=msbm5  scaled \magstep1
 \font\teneu=eufm10   scaled \magstep1
 \font\seveneu=eufm7 scaled \magstep1
 \font\fiveeu=eufm5 scaled \magstep1
 \font\tenib=cmmib10     scaled \magstep1
 \font\sevenib=cmmib7   scaled \magstep1
 \font\fiveib=cmmib5   scaled \magstep1
\fi
\newfam\msafam
\textfont\msafam=\tenmsa  \scriptfont\msafam=\sevenmsa
  \scriptscriptfont\msafam=\fivemsa
\newfam\msbfam
\textfont\msbfam=\tenmsb  \scriptfont\msbfam=\sevenmsb
  \scriptscriptfont\msbfam=\fivemsb
\def\Bbb{\ifmmode\let\next\Bbb@\else
 \def\next{\errmessage{Use \string\Bbb\space only in math mode}}\fi\next}
\def\Bbb@#1{{\Bbb@@{#1}}}
\def\Bbb@@#1{\fam\msbfam#1}
\newfam\eufam
\textfont\eufam=\teneu  \scriptfont\eufam=\seveneu
  \scriptscriptfont\eufam=\fiveeu
\def\frak{\ifmmode\let\next\frak@\else
 \def\next{\errmessage{Use \string\frak\space only in math mode}}\fi\next}
\def\frak@#1{{\frak@@{#1}}}
\def\frak@@#1{\fam\eufam#1}
\newfam\ibfam
\textfont\ibfam=\tenib  \scriptfont\ibfam=\sevenib
  \scriptscriptfont\ibfam=\fiveib
\def\bold{\ifmmode\let\next\bold@\else
 \def\next{\errmessage{Use \string\bold\space only in math mode}}\fi\next}
\def\bold@#1{{\bold@@{#1}}}
\def\bold@@#1{\fam\ibfam#1}
\def\hexnumber@#1{\ifcase#1 0\or 1\or 2\or 3\or 4\or 5\or 6\or 7\or 8\or
 9\or A\or B\or C\or D\or E\or F\fi}
\def\newsymbolb#1#2#3#4{\mathchardef#1="#2\hexnumber@\msbfam#3#4}

\relax

\hybrid

\begin{document}
\def\bea{\begin{eqnarray}}
\def\eea{\end{eqnarray}}
\def\beq{\begin{equation}}          \def\bn{\beq}
\def\eeq{\end{equation}}            \def\ed{\eeq}
\def\nn{\nonumber}                  \def\g{\gamma}
\def\Uq{U_q(\widehat{\frak{sl}}_2)}
\def\Uqp{U_q(\widehat{\frak{sl}}'_2)}
\def\Uqd{U^{*}_q(\widehat{\frak{sl}}_2)}
\def\uq{U_q({sl}_2)}
\def\uqd{U^*_q({sl}_2)}
\def\slaff{\frak{sl}^\prime_2}
\def\aff{\widehat{\frak{sl}}_2}
\def\ot{\otimes}
\def\id{\mbox{\rm id}}
\def\tr{\mbox{\rm tr}}
\def\th{\mbox{\rm th}}
\def\sh{\mbox{\rm sh}}
\def\ch{\mbox{\rm ch}}
\def\ctg{\mbox{\rm ctg}}
\def\tg{\mbox{\rm tg}}
\def\Re{{\rm Re}\,}
\def\RR{\Bbb{R}}
\def\ZZ{\Bbb{Z}}
\def\CC{\Bbb{C}}
\def\r#1{\mbox{(}\ref{#1}\mbox{)}}
\def\d{\delta}
\def\D{\Delta}
\def\da{{\partial_\alpha}}
\let\da=p
\def\R{{\cal R}}
\def\h{\hbar}
\def\Ga#1{\Gamma\left(#1\right)}
\def\ep{\varepsilon}
\def\ve{\ep}
\def\fract#1#2{{\mbox{\footnotesize $#1$}\over\mbox{\footnotesize $#2$}}}
\def\stackreb#1#2{\ \mathrel{\mathop{#1}\limits_{#2}}}
\def\res#1{\stackreb{\mbox{\rm res}}{#1}}
\def\lim#1{\stackreb{\mbox{\rm lim}}{#1}}
\def\Res#1{\stackreb{\mbox{\rm Res}}{#1}}
\let\dis=\displaystyle
\def\ee{{\rm e}}
\def\D{\Delta}
\renewcommand{\theequation}{{\thesection}.{\arabic{equation}}}
\def\Y-{\widehat{Y}^-}
\font\fraksect=eufm10 scaled 1728
\def\DYsect{\widehat{DY(\hbox{\fraksect sl}_2)}}
\def\DY{\widehat{DY(\frak{sl}_2)}}
\def\Yd{\DY}
\def\Ydd{\DY}
\let\z=\alpha
\let\b=\beta
\let\u=u
\let\v=v
\let\w=z
\def\g{\gamma}
\let\hsp=\qquad  
\def\trid{\tr\left(e^{\gamma d}\right)}

\begin{titlepage}
\begin{center}
\hfill ITEP-TH-8/96\\
\hfill BONN-TH-96-03\\
\hfill q-alg/9605039\\
\phantom.
\bigskip\bigskip\bigskip\bigskip\bigskip\bigskip
{\Large\bf Traces of Intertwining Operators\\
 for the Yangian Double}\\
\bigskip
\bigskip
{\large S. Khoroshkin\footnote{E-mail: khoroshkin@vitep1.itep.ru},
D. Lebedev\footnote{E-mail: lebedev@vitep1.itep.ru}}\\
\medskip
{\it Institute of Theoretical \& Experimental Physics\\
117259 Moscow, Russia}\\
\bigskip
{\large S. Pakuliak}\footnote{E-mail: pakuliak@thsun1.jinr.dubna.su}\\
\medskip
{\it Bogoliubov Laboratory of Theoretical Physics, JINR\\
141980 Dubna, Moscow region, Russia}\\
\medskip
{\it and}\\
\medskip
{\it Physikalisches Institut, Universit\"at Bonn\\
Nu\ss alle 12, 53115 Bonn, Germany}\\
\bigskip
\bigskip
\end{center}
\begin{abstract}
The traces over infinite dimensional representations of the central
extended Yangian double  for the product
of operators which intertwine these representations   are calculated.
For the special combinations of the intertwining operators the
traces are identified with form factors of local operators in 
$SU(2)$-invariant Thirring model. This identification is based
on the identities which are deformed analogs of the Gauss-Manin
connection identities for the hyperelliptic curves. 
\end{abstract}
\end{titlepage}
\clearpage
\newpage

\section{Introduction}

A big progress have been achieved last decades in the investigation of
 completely integrable quantum field theories.
By the complete integrability we understand the possibility to
write down explicit integral representation for an arbitrary 
form factor of any local operator in the model.
One of the first methods which allows one to gain complete integrability
was bootstrap program developed
for Sin-Gordon
and $SU(2)$-invariant Thirring models in the works by
F.A. Smirnov \cite{S1}. 
Each time when realization of the
bootstrap program led to success it always appeared as a miracle.
The reason for such a miracle lies in the deep
mathematical structure underlying the phyisical theories.
 One of  such mathematical
structures is a representation theory of infinite-dimensional Hopf 
algebras with the structure 
of doubles \cite{D}.

An application of such a symmetric ideology was realized
recently in the series of papers by Kyoto group \cite{XXZ,JMMN,JM}
devoted to investigation of 
quantum integrable lattice XXZ-type models in thermodinamic limit 
(for infinite lattices). The corresponding Hopf algebra for such
models have been identified with quantum deformation of universal
enveloping affine algebra $U_q(\widehat{\frak{g}})$ with $|q|<1$.
Since a representation theory for this Hopf algebras is studied
well enough, it is possible to calculate explicitly arbitrary physical
quantity in the model. The main ingredient of this approach was
a free field realization of infinite-dimensional integrable
representations of  $U_q(\widehat{\frak{g}})$ \cite{FJ}. 

The main goal of this paper is to expand this group-theoretical
ideology to the continuous integrable models. In principle, there  
exists a possibility to obtain the answers in continuous theory
using scaling limit from the corresponding 
lattice one \cite{L}. But it is possible to
do only for final objects, 
like integral formulas for form factors or correlation functions
of the local operators 
while some intermediate ones (representation
theory, for example) are lost in such a limit. 
We think that developing the direct group-theoretical 
methods in quantum integrable continuous field theories is an 
important problem.

Our main working example will be $SU(2)$-invariant Thirring model,
also known as  $SU(2)$ chiral Gross-Neveu model or current-current
perturbation of $SU(2)$ WZNW model at level 1. An infinite-dimensional
Hopf algebra associated with this model is a Yangian double \cite{B,BL,S2}.
Our main observation is that in order to develop group-theoretical
methods for  complete integration of $SU(2)$-invariant Thirring model
one has to use a central extension of the Yangian double 
$\widehat{DY(\frak{sl}_2)}$, introduced
recently in \cite{K} and investigated in \cite{KLP1}.
 For generalization to 
 $\frak{sl}_n$ and $\frak{gl}_n$ case see  \cite{I}. 

Organization of the paper is as follows. Section 2 is devoted to 
a short definition of central extension of the 
Yangian double $\DY$ and to realization of infinite dimensional
representations of it in terms of free bosonic field. In section 3
we describe the properties of type I and type II intertwining operators for 
these representations. Next section is devoted to definition
and calculation of traces for infinite-dimensional representation.
In sections 5 and 6 particular cases of the general trace formula
\r{general-trace} are considered. 
It is shown that traces of the products of type I
intertwining operators  coincide with 
correlation functions of the inhomogeneous XXX model \cite{N}. 
An analogous trace of type II operators is shown to be
equal identically zero
at ``physical'' value of shift parameter. 
In the last section the problem of
identification of group-theoretical and bootstrap
approaches to the calculation of form factors in completely integrable
field theories is considered.
This problem  was not cleared up even in case of lattice 
models associated with $U_q(\widehat{\frak{g}})$. 
We have found that 
this identification 
appears to be available because of special identities held between
Smirnov's integrals, which are deformed analogs of Gauss-Manin connection
identities for the  hyperelliptic integrals \cite{GM}.

\setcounter{equation}{0}            

\section{Central Extension of $\DYsect$}

Central extension of the double of
Yangian $DY(sl_2)$ is an infinite-dimensional  Hopf algebra 
  with central element $c$ and generators $d$,  
 $e_k, f_k, h_k$, $k\in {\ZZ}$,
gathered into generating functions \cite{K}
 \beq
\label{2.0}
e(\u)=\sum_{k\in\ZZ} e_k \u^{-k-1},
\hsp
f(\u)=\sum_{k\in\ZZ} f_k \u^{-k-1},
\hsp
h^{\pm}(\u)=1\pm\h\sum_{k\geq0\atop k<0}h_{k}\u^{-k-1},
\eeq
with the following relations
$$[d,e(\u)]=\frac{\mbox{d}}{\mbox{d}u}e(\u), \hsp
 [d,f(\u)]=\frac{\mbox{d}}{\mbox{d}u}f(\u), \hsp
 [d,h^\pm(\u)]=\frac{\mbox{d}}{\mbox{d}\u}h^\pm(\u),$$
\begin{eqnarray}
e(\u)e(\v)&=&{\u-\v+\h\over \u-\v-\h}\ e(\v)e(\u) \nn\\
f(\u)f(\v)&=&{\u-\v-\h\over \u-\v+\h}\ f(\v)f(\u) \nn\\
h^\pm(\u)e(\v)&=&{\u-\v+\h\over \u-\v-\h}\ e(\v)h^\pm(\u) \nn\\
h^+(\u)f(\v)&=&{\u-\v-\h-\h c\over \u-\v+\h-\h c}\ f(\v)h^+(\u) \nn\\
h^-(\u)f(\v)&=&{\u-\v-\h\over \u-\v+\h}\ f(\v)h^-(\u) \nn\\
h^+(\u)h^-(\v)
&=&{\u-\v+\h\over \u-\v-\h}\cdot{\u-\v-\h-\h c\over \u-\v+\h-\h c}\
h^-(\v)h^+(\u)
\nn\\
{[}e(\u),f(\v){]}&=&
{1\over \h}\left(\delta(\u-(\v+\h c))h^+(\u)-\delta(\u-\v)h^-(\v) \right)
\label{DY2}
\end{eqnarray}
where
$$\delta(\u-\v)=\sum_{n+m=-1}\u^n\v^m, \hsp
\delta(\u-\v)g(\u)=\delta(\u-\v)g(\v).$$

For the comultiplication formulas in $\widehat{DY(\frak{sl}_2)}$ see \cite{K},
 \cite{KLP1}.

Let ${\cal H}$ be a Heisenberg algebra generated by free
bosons with zero modes $a_{\pm n}$, $n=1,2,\ldots$,
$a_0$, $\da$  with commutation relations
$$
[a_n,a_m] =  n \delta_{n+m,0}\ ,\quad [\da,a_0]=2\ .
$$
 We use the following generating functions of elements from ${\cal H}$:
        \bea
        a_+(z)&=&\sum_{n\geq 1}\frac{a_n}{n}z^{-n}-p\log z\ ,\hsp
        a_-(z)=\sum_{n\geq 1}\frac{a_{-n}}{n}z^{n}+\frac{a_0}{2}\ ,
        \label{5.2} \\
        a(z)&=&a_+(z)-a_-(z),\hsp \phi_\pm(z)=\exp a_\pm(z)\ ,
        \label{5.2a}    \\
        {[}a_+(z),a_-(y){]}&=&-\log(z-y), \hsp |y|<|z|\ .
        \label{5.2b}
        \eea
        Let $V_i, i=0,1$ be formal power extensions of the Fock spaces
\beq
V_i={\CC}[[a_{-1},\ldots ,a_{-n},\ldots ]]\otimes
        \left(\oplus_{n\in{\ZZ}+i/2}{\CC}e^{na_0}\right)
\label{modules}
\eeq
with the action of bosons on these spaces
\bea
a_n&=&\hbox{the left multiplication by $a_n\otimes1$\ \  for $n<0$}\ ,\nn\\
   &=&{[}a_n,\ \cdot\ {]}\otimes 1\quad \hbox{for $n>0$}\ ,\nn\\
\ee^{n_1a_0}(a_{-j_k}\cdots a_{-j_1}\otimes\ee^{n_2a_0})&=&
a_{-j_k}\cdots a_{-j_1}\otimes\ee^{(n_1+n_2)a_0}\ ,\nn\\
u^{\da} (a_{-j_k}\cdots a_{-j_1}\otimes\ee^{na_0})&=&u^{2n}
a_{-j_k}\cdots a_{-j_1}\otimes\ee^{na_0}\ .\nn
\eea
\noindent
{\bf Proposition 1.}
{\it The following} $\mbox{End}\, V_i$-{\it 
valued functions satisfy commutation 
 relations \r{DY2} with $c=1$:}

\bea
e(\u)&=&\mu_-(\u-\h)\phi_+^{-1}(\u)\ ,\quad
 f(\u)=\mu_-^{-1}(\u)\phi_+(\u)\ ,\nn\\
h^+(\u)&=&\phi_+(\u-\h)\phi_+^{-1}(\u)\ ,\quad
 h^-(\u)=\phi_-(\u-\h)\phi_-^{-1}(\u+\h)\ ,
 \label{h-fields}                     \\
 e^{\g d}\phi_\pm(\u)&=&\phi_\pm(\u+\g)e^{\g d}\ ,\quad e^{\g d}(1\otimes 1)
 =1\otimes 1, \nn
\eea
where 
\beq
\mu_-(\u)=\phi_-(u+\h)\phi_-(u).
\label{q-kvad}
\eeq

Due to Proposition 1, the Fock spaces $V_i, i=0,1$ are $\DY$-modules, . 
We call  them basic (level $1$) representations of $\DY$.

\setcounter{equation}{0}

\section{Intertwining Operators for Basic 
 Representations }

Let $\eta_+(z)$ be the following $\mbox{End}\, V_i$-valued function (field):
\beq
\eta_+(z)=\lim{K\to\infty} (2\h K)^{-\da/2}
\prod_{k=0}^K {\phi_+(z-2k\h)\over\phi_+(z-\h-2k\h)}
\label{eta}
\eeq
The field $\eta_+(z)$ is well defined and satisfies a property 
\bea
\eta_+(z)\eta_+(z-\h)&=&\lim{K\to\infty}   (2\h K)^{-\da}
\prod_{k=0}^K {\phi_+(z-2k\h)\over\phi_+(z-2\h-2k\h)} \nn\\
&=&\phi_+(z) \lim{K\to\infty} (2\h K)^{-\da}
\phi^{-1}_+(z-2\h-2K\h) \nn                                 \\
&=&\phi_+(z) (-1)^{\da}
\label{prop-use}
\eea
where the operator $(-1)^{\da}$ acts on the module $V_i$ multipling it
by $(-1)^i$. We can rewrite $\eta_+(z)$ in an equivalent
form using the Stirling formula 
\bea
\eta_+(z)&=& (2\h)^{-\da/2}
\left( {\Ga{\fract{1}{2}-\fract{z}{2\h}}\over\Ga{-\fract{z}{2\h}}}    
\right) ^{-\da}
\nn\\
&\times&
\lim{K\to\infty}
\prod_{k=0}^{K}
\exp\left(\sum_{n=1}^{\infty}  {a_n\over n}
\left[(z-2k\h)^{-n}-(z-\h-2k\h)^{-n}\right]\right)                    
\label{eta1}
\eea
This presentation is useful to investigate the classical limit
of the field $\eta_+(z)$ 
\beq
\lim{\h\to 0}\eta_+(z)= \phi_+^{1/2}(z)(-1)^{\da/2}
\label{clas-eta}
\eeq
which shows that the formula \r{eta} can be treated as  a  
 deformed square root.

Let us define the intertwining
operators
\begin{eqnarray}
&\Phi^{(i)}(\z): V_i\to V_{1-i}\otimes V_\z,\quad
\Phi^{*(i)}(\z):V_i\otimes V_\z\to V_{1-i}
\ ,\label{typeI}\\
&\Psi^{(i)}(\b): V_i\to V_\b\otimes V_{1-i},\quad
\Psi^{*(i)}(\b):V_\b\otimes V_i\to V_{1-i}
\label{typeII}
\end{eqnarray}
which commute with the action of the Yangian double.
Here $V_z$ in \r{typeI} and \r{typeII} denotes two-dimensional
evaluation module
$$
V_z=V\otimes\CC[[z,z^{-1}]],\quad 
V=\CC v_+\otimes \CC v_-
$$
dependent on the spectral parameter $z$ \cite{KT}.
The components of the intertwinning operators are defined as follows
\bea
\Phi^{(i)}
(\z) v  &=&
\Phi^{(i)}_+(\z)v\otimes v_+ + \Phi^{(i)}_-(\z)v\otimes v_-\ ,\quad
\Phi^{*(i)}(\z)(v\otimes v_\pm)= \Phi^{*(i)}_\pm(\z)v \ ,
  \nn\\
\Psi^{(i)}(\b) v  &=&
v_+\otimes\Psi^{(i)}_+(\b)v + v_-\otimes\Psi^{(i)}_-(\b)v\ ,\quad
\Psi^{*(i)}(\b)(v_\pm\otimes v)= \Psi^{*(i)}_\pm(\b)v \ ,
 \nn  
\eea
where $v\in V_i$. Fix also the normalization of these operators
\bea\label{normaliz}
\Phi^{(i)}(\z)u_{i}&=&(-\z)^{-i/2}
\ep_{1-i}u_{{1-i}}
\otimes v_{\ep_i}+\cdots,\nn\\
\Psi^{*(i)}(\b)(v_{\ep_{1-i}}\ot u_{i})&=&(-\b)^{-i/2}u_{{1-i}}+\cdots,\quad
\ep_0=-,\ \ \ep_1=+\  ,
\eea
where we denote
$$u_0=1\ot 1\qquad \hbox{and} \qquad u_1=1\ot e^{a_0/2}\ . $$
Normalization condition \r{normaliz} allows us to write down
precise expressions for
$$
\Phi_\ep=\Phi^{(0)}_\ep\oplus\Phi^{(1)}_\ep:V_0\oplus V_1\rightarrow
V_1\oplus V_0
\quad\hbox{and}\quad
\Psi^*_\ep=\Psi^{*(0)}_\ep\oplus\Psi^{*(1)}_\ep:V_0\oplus V_1\rightarrow
V_1\oplus V_0
$$
without dependence on the index $i$. We have the following \cite{KLP1}
\medskip

\noindent
{\bf Proposition 2.} {\it Intertwinning operators \r{typeI}
 have the free field realiza\-ti\-on:}
\begin{eqnarray}
\Phi_-(\z)&=&
\phi_-(\z+\h)\eta_+^{-1}(\z)\ ,\label{4.5b}\\
\Phi_+(\z)&=&\Phi_-(\z)f_0-f_0\Phi_-(\z)\ ,
                \label{fI4} \\
   \Phi^{*(i)}_\nu(\z)&=&\nu(-1)^i\Phi^{(i)}_{-\nu}(\z-\h)\ ,\quad
   \nu=\pm\ ,
\label{comparI}\\
\Psi^{*}_-(\b)&=&
\phi_-^{-1}(\b)\eta_+(\b)\ ,\label{4.6b}\\
\Psi^{*}_+(\b)&=&e_0\Psi^*_-(\b)-\Psi^{*}_-(\b)e_0\ ,
              \label{f*II4}                                \\
   \Psi^{(i)}_\ep(\b)&=&\ep(-1)^{-i}\Psi^{*(i)}_{-\ep}(\b-\h) ,\quad
   \ep=\pm\ . \label{comparII}
\end{eqnarray}
\smallskip

These intertwining  operators satisfy commutation relations
\begin{eqnarray}
\Phi_{\ep_2}(\z_2) \Phi_{\ep_1}(\z_1)
&=&R_{\ep_1\ep_2}^{\ep'_1\ep'_2} (\z_1-\z_2)
\Phi_{\ep'_1}(\z_1) \Phi_{\ep'_2}(\z_2)\ , \label{ZFI} \\
\Psi^*_{\ep_1}(\b_1) \Psi^*_{\ep_2}(\b_2)
&=&-R_{\ep_1\ep_2}^{\ep'_1\ep'_2} (\b_1-\b_2)
\Psi^*_{\ep'_2}(\b_2) \Psi^*_{\ep'_1}(\b_1)\ , \label{ZFII} \\
\Phi_{\ep_1}(\z) \Psi^*_{\ep_2}(\b)
&=&\tau(\z-\b)
\Psi^*_{\ep_2}(\b) \Phi_{\ep_1}(\z)\ , \label{ZFI-II}\\
g\sum_{\ep}\Phi^{*(1-i)}_\ep(\z)\Phi^{(i)}_\ep(\z)&=&-\id
\ ,\label{inverse}\\
\label{orthI}
g\Phi^{(1-i)}_{\ep_1}(\z)\Phi^{*(i)}_{\ep_2}(\z)
&=&\delta_{\ep_1\ep_2}\ \id\ ,\\
\label{orthII}
g^{-1}\Psi^{(1-i)}_{\ep_1}(\b_1)\Psi^{*(i)}_{\ep_2}(\b_2)
&=&{\delta_{\ep_1\ep_2}
\over \b_1-\b_2} + o(\b_1-\b_2)\ ,
\end{eqnarray}
where  $R$-matrix is given by 
\beq\label{ver-R-mat}
R(z)=r(z)\overline R(z)
\eeq
and
\begin{eqnarray}
r(z)&=&{\Ga{\fract{1}{2}-\fract{z}{2\h}}\Ga{1+\fract{z}{2\h}}\over
        \Ga{\fract{1}{2}+\fract{z}{2\h}}\Ga{1-\fract{z}{2\h}}
   }\ ,\label{norm1}\\
\tau(z)&=&{
\Ga{1+\fract{z}{2\h}}\Ga{-\fract{z}{2\h}}\over
\Ga{\fract{1}{2}+\fract{z}{2\h}}
\Ga{\fract{1}{2}-\fract{z}{2\h}}}
\ =\ -\mbox{\rm ctg}\
{\pi z\over2\h}
   \ , \nn\\
\overline R(z)&=& \left(
\begin{array}{cccc}
1&0&0&0\\  0&b(z)&c(z)&0\\  0&c(z)&b(z)&0\\  0&0&0&1
\end{array}           \right)      \ ,        \label{R-mat}
\\
b(z)&=&{z\over z+\h},\quad c(z)\ =\  {\h\over z+\h}\ , 
\quad g\ =\ \sqrt{2\h\over\pi}\ . \nn
\end{eqnarray}
One can easily check that $R$-matrix \r{ver-R-mat}
satisfies the unitary and the crossing sym\-met\-ry conditions
$$
R(z)R(-z)=1 \ ,
$$
$$
(C\ot\id)\, R(z)\, (C\ot\id) =R^{t_1}(-z-h)
$$
with the charge conjugation matrix
\beq\label{charge}
C=\left(\begin{array}{cr} 0&-1\\1&0 \end{array}\right)
\eeq
and $R^{t_1}(z)$ means the transposition with respect to first space.

The verification of these commutation relations is based on
the following normal ordering relations (see details in \cite{KLP1})
\bea
\Phi_-(\z_2)\Phi_-(\z_1)      &=&
(2\h)^{1/2}{\Ga{1+\fract{\z_1-\z_2}{2\h}}\over
\Ga{\fract{1}{2}+\fract{\z_1-\z_2}{2\h}}}{:}\Phi_-(\z_2)\Phi_-(\z_1){:}\ ,
\label{fusionI}                        \\
\Psi^*_-(\b_1)\Psi^*_-(\b_2)      &=&
(2\h)^{1/2}{\Ga{\fract{1}{2}-\fract{\b_1-\b_2}{2\h}}\over
\Ga{-\fract{\b_1-\b_2}{2\h}}}{:}\Psi^*_-(\b_1)\Psi^*_-(\b_2){:}\ ,
\label{fusionII}                    \\
\Psi^*_\ep(\b)\Phi_\nu(\z)      &=&
(2\h)^{-1/2}{\Ga{\fract{1}{2}+\fract{\z-\b}{2\h}}\over
\Ga{1+\fract{\z-\b}{2\h}}}{:}\Psi^*_\ep(\b)\Phi_-(\z){:}\ ,
\label{fusionI-II}                    \\
\Phi_-(\z)f(\u)   &=&{1\over \u-\z}\    {:}\Phi_-(\z)f(\u){:}  \nn\\
f(\u)\Phi_-(\z)   &=&{1\over \u-\z-\h}\ {:}\Phi_-(\z)f(\u){:}  
\label{norderI}\\ 
\Psi^*_-(\b)e(\v) &=&{1\over \v-\b-\h}\ {:}\Psi^*_-(\b)e(\v){:}\nn\\ 
e(\v)\Psi^*_-(\b)  &=&{1\over \v-\b}\    {:}\Psi^*_-(\b)e(\v){:}
\label{norderII}\\ 
e(\v)\Phi_-(\z)  &=& \Phi_-(\z) e(\v) 
 = (\v-\z-\h)\ {:}\Phi_-(\z) e(\v){:}\nn\\ 
f(\u)\Psi^*_-(\b)&=& \Psi^*_-(\b)f(\u)= 
(\u-\b)\    {:} \Psi^*_-(\b)f(\u){:}\nn
\end{eqnarray}

These normal ordering formulas are used in particular to represent
the second components of the intertwinning operators where some
contour integrations over variables $\u$ and $\v$ are supposed.
Let us specify the contours in these integrals. By inspection
of \r{norderI} we can easily find that
the point $\z+\h$ should be inside the contour $C$
in the definition of the operator   $\Phi_+(\z)$
\beq
\Phi_+(\z)=
\Phi_-(\z)f_0-f_0\Phi_-(\z)
=\int_{C}{d\u\over 2i\pi}
\left(\Phi_-(\z)f(u)-f(u)\Phi_-(\z)\right)\ ,
\label{Phi+}
\eeq
while the point $\z$ should be  outside  $C$.
For the type II intertwinning
operators the picture is inverse.
The contour $\tilde C$ in definition of the operator
$\Psi^*_+(\b)$
\beq
\Psi^*_+(\b)
= e_0 \Psi^*_-(\b)-\Psi^*_-(\b) e_0
=\int_{\tilde C}{d\v\over 2i\pi}
\left(e(v)\Psi^*_-(\b)-\Psi^*_-(\b)e(v)\right)
\label{Psi+}
\eeq
should include the point $\b$ but exclude the point $\b+\h$.
   The form of the countours $C$ and
$\tilde C$ depends on the presence of other operators  and should
be specified in   each particular case of the product of 
intertwining operators.

\setcounter{equation}{0}

\section{Trace Formula}

The goal of this section is to obtain the trace of the product 
of arbitrary number of type I and type II intertwining operators.
We need the following 
\medskip

\noindent
{\bf Proposition 2.} 
{\it Let }
\beq
O=\prod_{j}\phi_-^{k_j}(w_j) \prod_{k}\eta_+^{p_k}(z_k) 
\label{O}
\eeq
{\it  
be an operator acting in the infinite dimensional Fock space $V_0$.
Under the restrictions for the integers 
$k_j$ and $p_k$ }
\beq
\sum_j k_j=0 \qquad  \sum_k p_k=0
\label{cond-gener}
\eeq
{\it a ratio} 
$${\tr\left(e^{\g d}O\right)\over
\tr\left(e^{\g d}\right)}$$ {\it
of traces over this infinite dimensional
representation space
   can be calculated
using the formula}
\beq
{\tr\left(e^{\gamma d}O\right)\over
\tr\left(e^{\g d}\right)}
= \prod_{m=1}^\infty
\exp\left(\res{x,y}{{<}a_+(x-m\g)[O-1]a'_-(y) {>}\over x-y} \right)
\label{trace}
\eeq
\smallskip
{\it and is equal to}
\beq
{1\over \tr\left(e^{\gamma d}\right)}
\tr\left(e^{\gamma d}
\prod_j\phi_-^{k_j}(w_j) \prod_k\eta_+^{p_k}(z_k) \right)=
\prod_{m=1}^\infty \prod_{k=0}^\infty \prod_{j,k}
{(z_k-w_j-\h-m\g-2k\h)^{k_jp_k} \over
(z_k-w_j-m\g-2k\h)^{k_jp_k}}\ ,
\label{calc-gener}
\eeq
{\it where an integral representation for the r.h.s. of} \r{calc-gener}
{\it is} $(\mbox{Re}\,\g/\h>0) $  
\beq 
\exp\left({1\over4}
\int_{0}^\infty{dx\over x}
\sum_{j,k}{k_jp_k}
{e^{\fract{(\h-\g)x}{2\h}}e^{\fract{(z_k-w_j)x}{\h}}
\over\sh\fract{\g x}{2\h}\ch\fract{x}{2}
}\right) .
\label{int-rep}
\eeq

Formula \r{trace} is an analog of the formula (8.6) given in the book
\cite{JM} for $U_q(\widehat{\frak{sl}}_2)$. To obtain \r{trace}
one has to pass from summation over powers of different modes of free
bosonic field, to
the summation over products of the free field taken in different points.
In the basis constructed from these fields an operator
$e^{\gamma d}$ acts by shift of the spectral parameter. 
The proof of \r{trace} is complicated but straightforward and
can be found 
in \cite{Ch}. A significant difference of trace calculations
in Yangian case in contrast to quantum affine one \cite{JM} is that
only ratio of traces is well defined in a sense that 
corresponding infinite products are convergent.

Because of the property \r{prop-use} of the field $\eta_+(z)$ 
and definition \r{q-kvad} of the field $\eta_-(z)$ 
the formula \r{calc-gener} is sufficient to calculate 
\r{trace} when 
the operator $O$ is a product of different combinations
of free fields $\phi_-(z)$,
$\phi_+(z)$,
$\mu_-(z)$ and $\eta_+(z)$.

Let us apply above formulas for computation of
 the trace over the space $V_0$
of the following product of type I and type II intertwining operators
\beq
\Psi^*_{\ep_N}(\b_N)\ldots\Psi^*_{\ep_1}(\b_1)
\Phi_{\nu_M}(\z_M)\ldots\Phi_{\nu_1}(\z_1)
\label{prod}
\eeq
with the numbers of ``plus'' components of the operators 
$$
\{\#\,j\,|\,\ep_j=+\}=n\qquad \{\#\,i\,|\,\nu_i=+\}=m\ .
$$
First, we have to normal order the product \r{prod} using formulas of the
previous section. We can rewrite \ref{prod}  ($\z_{ij}=\z_i-\z_j$ 
and $\b_{ij}=\b_i-\b_j$ ) as following:
\bea 
&\prod_{j=1}^N\Psi^*_{\ep_j}(\b_j)\prod_{i=1}^M\Phi_{\nu_i}(\z_i)=
(-\h)^{n+m}(2\h)^{{N(N-1)+M(M-1)-2NM\over 4}}\nn\\
& \times \prod_{i<j}
{ \Ga{\fract{1}{2}+\fract{\b_{ij}}{2\h}} \over
\Ga{\fract{\b_{ij}}{2\h}} } 
\prod_{i<j}
{ \Ga{1+\fract{\z_{ij}}{2\h}} \over
\Ga{\fract{1}{2}+\fract{\z_{ij}}{2\h}} } 
\prod_{j=1}^N\prod_{i=1}^M
{ \Ga{\fract{1}{2}+\fract{\z_i-\b_j}{2\h}} \over
\Ga{1+\fract{\z_i-\b_j}{2\h}} } \nn\\
&\times
\stackreb{\int\!\!\cdots\!\!\int}{\tilde C}
{d\v_k\over2\pi i}
{\tilde P_{\ep}(\v;\b)\prod_{i<j}(\v_i-\v_j)(\v_i-\v_j-\h)\over
\prod_{i=1}^n\prod_{j=1}^N (\v_i-\b_j)(\v_i-\b_j-\h)}\nn\\
&\times
\stackreb{\int\!\!\cdots\!\!\int}{C}
{d\u_p\over2\pi i}
{P_{\nu}(\u;\z)\prod_{i<j}(\u_i-\u_j)(\u_i-\u_j+\h)\over
\prod_{i=1}^m\prod_{j=1}^M (\u_i-\z_j)(\u_i-\z_j-\h)}\nn\\
&\times
{\prod_{k=1}^n \prod_{i=1}^M (\v_k-\z_p-\h)
\prod_{p=1}^m \prod_{j=1}^N (\u_p-\b_j)
\over
\prod_{k=1}^n\prod_{p=1}^m 
 (\u_p-\v_k)(\u_p-\v_k+\h)}\ O_-O_+\ ,\nn
\eea  
where
\bea
O_-&=&
\prod_{j=1}^N\phi_-^{-1}(\b_j)\prod_{i=1}^M\phi_-(\z_i+\h)
\prod_{k=1}^n\mu_-(\v_k-\h)   \prod_{p=1}^m\mu_-^{-1}(\u_p)
\label{O-}\\
O_+&=&
\prod_{j=1}^N\eta_+(\b_j)\prod_{i=1}^M\eta_+^{-1}(\z_i)
\prod_{k=1}^n\phi_+^{-1}(\v_k)\prod_{p=1}^m\phi_+(\u_p)
\label{O+}
\eea
and for 
\beq
\{b_1<\ldots<b_n\}=\{j\,|\,\ep_j=+\}
\quad {\rm and}\quad
\{a_1<\ldots<a_n\}=\{i\,|\,\nu_i=+\}
\label{index}
\eeq
we define the polynomials
\bea
P_{\nu}(\u;\z)&=&\prod_{p=1}
\prod_{i<a_p}(\u_p-\z_i)
\prod_{i>a_p} (\u_p-\z_i-\h)\ ,\nn\\
\tilde P_{\ep}(\v;\b)&=&\prod_{k=1}
\prod_{j>b_k}(\v_k-\b_j)
\prod_{j<b_k} (\v_k-\b_j-\h)\ .\nn
\eea

Applying formula \r{calc-gener} to the operators
$O_-O_+$ given by \r{O-} and \r{O+} we see
that the condition \r{cond-gener}
for the existence of this trace 
can be written as
\beq
{N-M\over2}=n-m
\label{condition}
\eeq 
and the trace of the product \r{prod} is equal to
\bea 
&{1\over \tr\left(e^{\g d}\right)}
\tr\left(e^{\g d}
\prod_{j=1}^N\Psi^*_{\ep_j}(\b_j)\prod_{i=1}^M\Phi_{\nu_i}(\z_i)
\right)\nn\\
&\quad=A_{NM}^{nm}\prod_{j<j'}\tilde\zeta(\b_j-\b_{j'})
\prod_{i<i'}\zeta(\z_i-\z_{i'})
\prod_{i,j}\overline\zeta(\z_i-\b_{j})\nn\\
&\qquad \times \prod_{k=1}^n \int_{\tilde C}{d\v_k\over2\pi i}
\prod_{j=1}^{N}
\Ga{{\b_j-\v_i\over\g}}\Ga{{\v_i-\b_j-\h\over\g}}
{\tilde P_{\ep}(\v;\b)\over \g^{n(N-1)}}\nn\\
&\quad\qquad \times
\prod_{k<k'}
{\sin\pi\fract{\v_k-\v_{k'}}{\g}\over\Ga{\fract{\v_k-\v_{k'}-\h}{\g}}
\Ga{\fract{\v_{k'}-\v_k-\h+\g}{\g}}}\nn\\
&\qquad \times \prod_{p=1}^m \int_{C}{d\u_p\over2\pi i}
\prod_{i=1}^{M}
\Ga{{\u_p-\z_i\over\g}}\Ga{{\z_i-\u_p+\h\over\g}}
{P_{\nu}(\u;\z)\over \g^{m(M-1)}} \nn\\
&\quad\qquad \times
\prod_{p<p'}
{\sin\pi\fract{\u_p-\u_{p'}}{\g}\over\Ga{\fract{\u_p-\u_{p'}+\h}{\g}}
\Ga{\fract{\u_{p'}-\u_p+\h+\g}{\g}}}\nn\\
&\quad\qquad
\times 
{
\prod_{p=1}^m\prod_{j=1}^N\sin\pi\fract{\u_p-\b_j}{\g}
\prod_{k=1}^n\prod_{i=1}^M\sin\pi\fract{\v_k-\z_i}{\g}
\over
\prod_{p=1}^m\prod_{k=1}^n\sin\pi\fract{\u_p-\v_k}{\g}
\sin\pi\fract{\u_p-\v_k+\h}{\g}     }
\label{general-trace}
\eea
This formula in particular case $\g=2\h$ has been written in \cite{L}
and can be obtained by means of the scaling limit from
general trace formula given in \cite{JM} for the quantum affine algebra.

In \r{general-trace} the numerical constant $A_{NM}^{nm}$ is
\beq
A_{NM}^{nm}= \left({2\h\over\sqrt\pi\g}\right)^{(N-M)^2/4}
{(-1)^{Nn+Mm}(-\h)^{n+m}\pi^{nm+n/2+m/2}\tilde G^{N/4}(0) G^{M/4}(0)
\over(2\h)^{(N+M)/4}\g^{2n+2m}
\Ga{1-\fract{\h}{\g}}^{(N^2+4n)/4}\Ga{1+\fract{\h}{\g}}^{(M^2+4m)/4}
}
\label{const-gen}
\eeq
and we introduced functions
\bea
\zeta(\z)&=& \left[
{ \Ga{1+\fract{\z}{2\h}} \over
\Ga{\fract{1}{2}+\fract{\z}{2\h}} }
{G(\z)\over G^{1/2}(0)}
\right]\ ,\nn\\
\tilde \zeta(\b)&=& \left[
{ \Ga{\fract{1}{2}+\fract{\b}{2\h}} \over
\Ga{\fract{\b}{2\h}} }
{\tilde G(\b)\over \tilde G^{1/2}(0)}
\right]\ ,\nn\\
\overline\zeta(z)&=& \left[
{ \Ga{\fract{1}{2}+\fract{z}{2\h}} \over
\Ga{1+\fract{z}{2\h}} }
\overline G(z)
\right]\nn
\eea
and 
\bea
G(\z)&=&
\prod_{m=1}^\infty\prod_{k=0}^\infty
{(-\h-2k\h-m\g) (-3\h-2k\h-m\g) \over
(-2k\h-m\g) (-2\h-2k\h-m\g)  }\nn\\
&&\times
{(\z-\h-2k\h-m\g) (-\z-\h-2k\h-m\g)
\over (\z-2\h-2k\h-m\g) (-\z-2\h-2k\h-m\g) }\ ,
\label{G-funI}\\
\tilde G(\b)&=&
\prod_{m=1}^\infty\prod_{k=0}^\infty
{(-2k\h-m\g) (-2\h-2k\h-m\g) \over
(-\h-2k\h-m\g) (\h-2k\h-m\g)  }\nn\\
&&\times
{(\b-2k\h-m\g) (-\b-2k\h-m\g)
\over (\b-\h-2k\h-m\g) (-\b-\h-2k\h-m\g) }\ ,
\label{G-funII}\\
\overline G(z)&=&
\prod_{m=1}^\infty\prod_{k=0}^\infty
{(-2k\h-m\g) \over  (-2\h-2k\h-m\g)   }\nn\\
&&\times{(z-\h-2k\h-m\g) (-z-2\h-2k\h-m\g)
\over (z-2k\h-m\g) (-z-\h-2k\h-m\g) }\ .
\label{G-funI-II}
\eea

Functions \r{G-funI}, \r{G-funII} and \r{G-funI-II}
are well defined because the infinite product
of rational functions \cite{Ba}
\beq
\prod_{k_1,\ldots,k_n\geq0}
{\prod_m(a_m+\sum k_j \omega_j)\over\prod_p(b_p+\sum k_j \omega_j)}
\label{Barnes}
\eeq
is well defined as function of the variables $a_m$ and $b_p$
if the following constraints for these variables are satisfied
\beq
\sum_m (a_m)^q=\sum_p (b_p)^q,\qquad q=0,1,\ldots,n\ .
\label{constr}
\eeq
In particular, for $n=1$ the function \r{Barnes} is equal to the
ratio of $\Gamma$-functions. The simplest way to prove \r{constr} is to
use an integral representation for the function \r{Barnes}.

 The functions 
$G(\z)$,  $\tilde G(\b)$ and  $\overline G(z)$  
have the following integral representations in appropriate regions of 
 parameters $\alpha, \beta , z$ 
 (Re$\,\g/\h>0$)
\bea
G(\z)&=&
\exp\left(\int_0^\infty{dx\over x}
{\th\fract{x}{2}\sh\fract{(\z+\h)x}{2\h}\sh\fract{(\z-\h)x}{2\h}
\over
\sh\fract{\g x}{2\h}}
\exp\left({-\frac{(\g-2\h)x}{2\h}}\right)\ ,
\right)\label{iI}\\
\tilde G(\b)&=&
\exp\left(\int_0^\infty{dx\over x}
{\th\fract{x}{2}\sh\fract{(\b+\h)x}{2\h}\sh\fract{(\b-\h)x}{2\h}
\over
\sh\fract{\g x}{2\h}}
\exp\left({-\frac{\g x}{2\h}}\right)\ ,
\right)\label{iII}\\
\overline G(z)&=&
\exp\left(-\int_0^\infty{dx\over x}
{\th\fract{x}{2}\sh\fract{(z+\h)x}{2\h}\sh\fract{z x}{2\h}
\over
\sh\fract{\g x}{2\h}}
\exp\left({\fract{(\g-\h)x}{2\h}}\right)
\right)
\label{iI-II}
\eea
and they satisfy the  properties
\beq
G(z)=G(-z),\quad G(\h)=1,\quad \tilde G(z)=\tilde G(-z),\quad 
\tilde G(\h)=1,\quad \overline G(0)=1,\label{p1}
\eeq
\beq
G(z)\tilde G(z)={z\over\h}{\sin(\pi\h/\g)\over \sin(\pi z/\g)},\quad
\overline G(z)\overline G(-z)= \overline G(z)\overline G(z-\h)= 
{\pi z\over\g}{1\over\sin(\pi z/\g)}\ .
\label{p2}
\eeq
Properties \r{p1} are useful for checking accordance of the 
formula \r{general-trace} with the normalization of the 
intertwining operators
\r{inverse}--\r{orthII}. 

The last thing which we should to do is to specify the contours
$C$ and $\tilde C$ in \r{general-trace}. The correct choice of the 
contours is dictated by \r{Phi+}, \r{Psi+} and prescriptions how
to use the formula \r{calc-gener}. 
The contour $C$ for the 
 integration over $u_p$, $p=1,\ldots,m$ is such that 
($r=0,1,2,\ldots$)
\bea
&\z_j+\h+r\g\quad\hbox{are inside  $C$},\quad
\z_j-r\g\quad\hbox{are outside  $C$},\quad 
 j=1,\ldots,M\nn\\
&\v_k-\h+\g+r\g,\quad\v_k+\g+r\g
 \quad\hbox{are inside  $C$},\nn\\
&\v_k-\h-r\g,\quad \v_k-r\g
\quad\hbox{are outside  $C$},\quad 
 k=1,\ldots,n  \label{conI}
\eea
and the contour $\tilde C$ for the 
 integration over $v_k$, $k=1,\ldots,n$ 
is such that ($r=0,1,2,\ldots$)
\bea
&\b_i+r\g\quad\hbox{are inside  $\tilde C$},\quad
\b_i+\h-r\g\quad\hbox{are outside  $\tilde C$},\quad 
 i=1,\ldots,N\nn\\
&\u_p+\h+r\g,\quad\u_p+r\g
 \quad\hbox{are inside  $\tilde C$},\nn\\
&\u_p+\h-\g-r\g,\quad \u_p-\g-r\g
\quad\hbox{are outside  $\tilde C$},\quad 
 p=1,\ldots,m\ .  \label{conII}
\eea
Also we should remark  that both contours
cross the infinity  point along the line
$z=i\h t$, $t\in\RR$ in the complex planes corresponding 
to the integration variables $v_k$ and $u_p$
because they are pinched 
from both sides of this line 
by the poles of the $\Gamma$-functions.
If the integrand in the general formula \r{general-trace}
has a pole in the point $\infty$ then we should understand
corresponding integrals as principal value integrals
(see next section).

It is quite difficult to work with general formula \r{general-trace}
so in next four sections  we will consider  particular cases
of this formula which are of special interest for us.

\setcounter{equation}{0}

\section{Correlation Functions}

The experience of group-theoretical approach in quantum integrable
lattice models \cite{JM} teaches 
that some combinations of type I intertwining 
operators can be identified with generating functions of the 
local operators in the corresponding 
quantum integrable model. 
Therefore the calculation of traces of type I intertwining operators gives 
 in particular  
the  correlation functions 
of local operators.
In this  case we have
from \r{general-trace} for $N=n=0$ and $M=2m$ due to \r{condition}
\bea
&{1\over\trid}
\tr\left(e^{\g d}
\prod_{i=1}^{2m} \Phi_{\nu_i}(\z_i)
\right) = A_{0,2m}^{0,m}
\prod_{i<i'} \zeta(\z_i-\z_{i'}) 
\nn\\
&\qquad \times \prod_{p=1}^m \int_{C}{d\u_p\over2\pi i}
\prod_{i=1}^{2m}
\Ga{{\u_p-\z_i\over\g}}\Ga{{\z_i-\u_p+\h\over\g}}
{P_{\nu}(\u;\z)\over \g^{m(2m-1)}} \nn \\
&\qquad \times
\prod_{p<p'}
{\sin\pi\fract{\u_p-\u_{p'}}{\g}\over\Ga{\fract{\u_p-\u_{p'}+\h}{\g}}
\Ga{\fract{\u_{p'}-\u_p+\h+\g}{\g}}}\label{Nakaya}
\eea

For $\g=2\h$ the  formula \r{Nakaya} was  obtained in \cite{N}
by means of scaling limit from the corresponding formula in XXZ model.
In this particular case there exist a possibility to 
reduce the number of integrals in \r{Nakaya} by a method suggested 
by F.A. Smirnov \cite{N}

We  demonstrate this in the particular case of $m=1$. Integral 
in \r{Nakaya} is interesting  in this simplest case   because
it can be calculated explicitely for arbitrary step $\g$.
We have in this particular case instead of \r{Nakaya} the following
integral ($\z=\z_1-\z_2$)
\bea
&{1\over\trid}
\tr\left(e^{\g d}\Phi_{\nu_2}(\z_2)\Phi_{\nu_1}(\z_1)\right)=
{(2\h)^{1/2}   
(-\h)\over \g^3 \Gamma^2
\left(\fract{\g+\h}{\g}\right)}
{\Ga{1+\fract{\z_1-\z_2}{2\h}}\over
\Ga{\fract{1}{2}+\fract{\z_1-\z_2}{2\h}}}G(\z) 
\nn\\
&\qquad \times \int_{C}  {d\u\over 2\pi i}
{P_{\nu_2\nu_1}(\u)\over\g}
\prod_{j=1}^2  \Ga{\frac{\u-\z_j}{\g}} \Ga{\frac{\z_j+\h-\u}{\g}}
\label{tr2I}
\eea

Denoting the  integral in \r{tr2I} by $I_{\nu_2\nu_1}$ we can obtain 
that 
\beq
I_{\nu_2\nu_1} =  -\nu_2 \g
{\Ga{\fract{\g+\h}{\g}} \Ga{\fract{\h}{\g}}
\Ga{\fract{\z+\h}{\g}} \Ga{\fract{\g+\h-\z}{\g}} \over
\Ga{\fract{\g+2\h}{\g}}          }
\label{int2I}
\eeq
because of the following classical result 
(Mellin-Barnes type integral)  \cite{BE}
\beq
\int_{-i\infty}^{+i\infty} {ds\over 2\pi i}
\Ga{a+s}\Ga{b+s}\Ga{c-s}\Ga{d-s}\nn\\
={\Ga{a+c}\Ga{a+d}\Ga{b+c}\Ga{b+d}
\over\Ga{1+a+b+c+d}}
\label{classical}
\eeq
Hence
\beq
{1\over\trid}
\tr\left(e^{\g d}\Phi_{\nu_2}(\z_2)\Phi_{\nu_1}(\z_1)\right)=
{\nu_2 G(\z)\over(2\h)^{1/2}}   {\Ga{1+\fract{\z}{2\h}}\over
\Ga{\fract{3}{2}+\fract{\z}{2\h}}}
{\Ga{\fract{\g+\h+\z}{\g}} \Ga{\fract{\g+\h-\z}{\g}} \over
\Ga{\fract{\g+2\h}{\g}}          }\ .  \label{fin2I}
\eeq

Applicability of the formula \r{classical} to calculate  the
integral $I_{\nu_2\nu_1}$
should be checked when $\g =2\h$, because the
integrand has a pole at the point $\infty$ in this case.
Setting $\g=2\h$ in \r{fin2I} we obtain
\beq
{1\over\tr\left(e^{2\h d}\right)}
\tr\left(e^{2\h d}
\Phi_{\nu_2}(\z_2)\Phi_{\nu_1}(\z_1)\right) =
\nu_2 (2\h)^{-1/2}   G(\z)
\Ga{1+\frac{\z}{2\h}} \Ga{\frac{3}{2}-\fract{\z}{2\h}}\ .
\label{int2Ilim}
\eeq
On the other hand we can calculate the
 integral $I_{\nu_2\nu_1}\bigr|_{\g=2\h} $
using Smirnov's hint. 
Using the identity
\bea
\cos\pi\frac{\z_1-\z_2}{2\h} &=&
\cos\pi\frac{\u-\z_1}{2\h} \cos\pi\frac{\u-\z_2}{2\h} +
\sin\pi\frac{\u-\z_1}{2\h} \sin\pi\frac{\u-\z_2}{2\h} \nn\\
&=& \sin\pi\frac{\u-\z_1-\h}{2\h} \sin\pi\frac{\u-\z_2-\h}{2\h} +
\sin\pi\frac{\u-\z_1}{2\h} \sin\pi\frac{\u-\z_2}{2\h}
\label{Smir-ide}
\eea
and the formula
$$\Ga{x}\Ga{1-x}={\pi\over\sin\pi x}$$
we can reduce this integral as follows
\bea
I_{\nu_2\nu_1}&=&
{\pi^2\over\cos\pi\fract{\z_1-\z_2}{2\h}}
\int_{C}  {d\u\over2\pi i} {P_{\nu_2\nu_1}(u)\over2\h}
\prod_{j=1}^2
{\Ga{\fract{\h+\z_j-\u}{2\h}}\over\Ga{\fract{2\h+\z_j-\u}{2\h}}}
\label{non-zero}\\
&&+\
{\pi^2\over\cos\pi\fract{\z_1-\z_2}{2\h}}
\int_{C}  {d\u\over2\pi i} {P_{\nu_2\nu_1}(u)\over2\h}
\prod_{j=1}^2 {\Ga{\fract{\u-\z_j}{2\h}}\over\Ga{\fract{\h+\u-\z_j}{2\h}}}
\label{zero}
\eea
Let us calculate first $I_{+-}$ since arguments for the second
integral $I_{-+}$ are the same. Using explicit form of the polynomial 
$P_{+-}(u)$ we can find that the integral in question
is equal to the sum of two integrals 
\bea
I_{{+}{-}}
&=& {\pi^2\over\cos\pi\fract{\z_1-\z_2}{2\h}}
\int_{C_2}{d\u\over 2\pi i}
{\Ga{\fract{1}{2}+\fract{\z_1-\u}{2\h}}\over\Ga{\fract{\z_1-\u}{2\h}}}
{\Ga{\fract{\h+\z_2-\u}{2\h}}\over\Ga{\fract{1}{2}+\fract{\h+\z_2-\u}{2\h}}}
\label{1}\\
&+&
{\pi^2\over\cos\pi\fract{\z_1-\z_2}{2\h}}
\int_{C_1}{d\u\over 2\pi i}
{\Ga{\fract{1}{2}+\fract{\h+\u-\z_1}{2\h}}\over
\Ga{\fract{\h+\u-\z_1}{2\h}}}
{\Ga{\fract{\u-\z_2}{2\h}}\over\Ga{\fract{1}{2}+\fract{\u-\z_2}{2\h}}}
\label{2}\\
&=& -{\pi^2\over\cos\pi\fract{\z}{2\h}}
{\z-\h\over2}=\pi\h
\Ga{\frac{1}{2}+\frac{\z}{2\h}} \Ga{\frac{3}{2}-\frac{\z}{2\h}}\nn
\eea
where contours $C_1$ and $C_2$ being small semicircles around the
point $\infty$ are shown on the Fig. 1.
\vskip0.5cm
\setlength{\unitlength}{0.00083333in}
\begin{picture}(4478,1185)(-1200,-10)
\thicklines
\path(189.000,531.000)(309.000,561.000)(189.000,631.000)
\put(309.000,861.000){\arc{600.000}{1.5708}{4.7124}}
\path(4029.000,1191.000)(3909.000,1161.000)(4029.000,1091.000)
\put(3909.000,861.000){\arc{600.000}{4.7124}{7.8540}}
\put(50,250){\makebox(0,0)[lb]{$C_1$}}
\put(309,800){\makebox(0,0)[lb]{$\cdot$}}
\put(3909,800){\makebox(0,0)[lb]{$\cdot$}}
\put(450,800){\makebox(0,0)[lb]{$\infty$}}
\put(3650,800){\makebox(0,0)[lb]{$\infty$}}
\put(3850,250){\makebox(0,0)[lb]{$C_2$}}
\put(1900,0){\makebox(0,0)[lb]{Fig.~1.}}
\end{picture}
\vskip0.7cm
Since the  arguments in ratios of $\Gamma$-functions \r{1}, \r{2}
differ by $1/2$ we can calculate these integrals 
using the Stirling formula. The same calculation for the integral $I_{-+}$
shows that formula  \r{fin2I}  is correct for all
possible values of $\g$.

Also one can easily check that \r{fin2I} is in accordance with 
the normalization
conditions for the type I intertwining operators \r{inverse} and 
\r{orthI} which 
in this case looks as
$$\Phi_{\nu_2}(\z)\Phi_{\nu_1}(\z-\h)=\nu_2 g^{-1} \delta_{\nu_2,-\nu_1}$$
and
$$   \sum_\nu \nu \Phi_{\nu}(\z-\h) \Phi_{-\nu}(\z)=g^{-1}\ . $$

\setcounter{equation}{0}

\section{Calculation of Multi-Point Traces for Type II VO}

The goal of this section is to prove the statement that in the 
``physical''
limit $\g=2\h$ the trace of the product of type II intertwining
operators is equal identically zero.
In such a case we have
from \r{general-trace} for $M=m=0$ and $N=2n$ due to \r{condition}:
\bea
&{1\over
\tr\left(e^{2\h d}\right)} 
 \tr\left(e^{2\h d}
\prod_{j=1}^{2n} \Psi^{*}_{\ep_i}(\b_i)
\right)
=   A_{2n,0}^{n,0} 
\prod_{j<j'} \tilde\zeta(\b_j-\b_{j'})
\nn\\
&\qquad \times \prod_{k=1}^n \int_{\tilde C}{d\v_k\over2\pi i}
\prod_{j=1}^{2n}
\Ga{{\b_j-\v_k\over2\h}}\Ga{{\v_k-\b_j-\h\over2\h}}
{\tilde P_{\ep}(\v;\b)\over (2\h)^{n(2n-1)}} \nn \\
&\qquad \times
\prod_{k<k'}
{\sin\pi\fract{\v_k-\v_{k'}}{2\h}\over
\Ga{\fract{\v_k-\v_{k'}-\h}{2\h}}
\Ga{\fract{\v_{k'}-\v_k+\h}{2\h}}}\label{IImnogo}
\eea
Using the identity
\bea
&\prod_{i<j}
{\sin\pi{\v_i-\v_j\over2\h}\over\Ga{{\v_i-\v_j-\h\over 2\h}}
\Ga{{\v_j-\v_i+\h\over2\h}}}= \left(\pi i\over 4\h \right)^{n(n-1)/2}
{2^{2n-1}\pi^{2n}\over\tilde B_1 B^{1/2}_{2n}}\nn\\
&\quad \times
\sum_{k=1}^{n} \Delta_k(v_j;j\neq k) \prod_{i<j\atop i,j\neq k}
(\v_i-\v_j-\h)
\prod_{j=1}^{k-1} (\v_k-\v_j+\h)
\prod_{j=k+1}^{n} (\v_k-\v_j-\h) \nn\\
&\quad \times
\left[
{1\over \prod_{j=1}^{2n}
\Ga{\fract{\v_k-\b_j-\h}{2\h}} \Ga{\fract{\b_j-\v_k+3\h}{2\h}} }
-{(-1)^n\over  \prod_{j=1}^{2n}
\Ga{\fract{\v_k-\b_j}{2\h}} \Ga{\fract{\b_j-\v_k+2\h}{2\h}} }
\right]
\nn
\eea
where
\bea
\tilde B_1&=& \sum_j e^{-\pi i\b_j/\h},\qquad
B_{2n}=\exp \left({\pi i\over\h}\sum_j\b_j\right)\nn\\
\Delta_k(v_j;j\neq k)&=&\det\left(e^{-(n-2l+1)\pi iv_j/\h}\right)
_{2\leq l\leq n\ \ 1\leq j\neq k\leq n}\nn
\eea
we can rewrite the integral in \r{IImnogo} over the variable $v_k$
as a difference of two integrals
\beq
\int_{\tilde C}{d\v_k\over 2\pi i}  f(\v_k)\left[
\prod _{j=1}^{2n}{\Ga{\fract{\b_j-\v_k}{2\h}}
\over\Ga{\fract{3\h+\b_j-\v_k}{2\h}}}
-(-1)^n
\prod _{j=1}^{2n}{\Ga{\fract{\v_k-\b_j-\h}{2\h}}
\over\Ga{\fract{\v_k-\b_j+2\h}{2\h}}}     \right]
\label{dif-int}
\eeq
where we have introduced a polynomial $f(\v_k)$
$$
f(\v)=
\prod_{j=1}^{k-1}(\v_k-\b_j-\h)
\prod_{j=k+1}^{2n}(\v_k-\b_j)
\prod_{j=1}^{k-1}(\v_k-\v_j+\h)
\prod_{j=k+1}^{n}(\v_k-\v_j-\h)
$$
We can see now that the integral in \r{dif-int} 
is equal to the sum of two integrals over two small semicircles around
the infinity point. These integrals are reduced to residues of the integrand
at this point. 
Using the Stirling formula one can easily
notice that the expansion of the integrand in \r{dif-int} in the vicinity
of the point $\infty$ starts from $v_k^{-2}$. Indeed,
the total power of polynomial
$f(\v_k)$ is equal $3n-2$, while the ratio of $\Gamma$-functions
produce $\v_k^{3n}$ in the denominator. Since it is true
for all terms in the sum over $k$ we proved that the trace
\beq
{1\over\tr\left(e^{2\h d}\right)}
\tr\left(e^{2\h d}\prod_{j=1}^{2n} \Psi^*_{\ep_i}(\b_j)
\right)     =0
\label{trII=zero}
\eeq

It is instructive 
to consider again the simplest case of the formula \r{IImnogo}
for $n=1$. 
This integral can be calculated using \r{classical}
and is equal to ($\b=\b_1-\b_2$)
\beq
{1\over\tr\left(e^{2\h d}\right)}
\tr\left(e^{\g d}\Psi^*_{\ep_2}(\b_2)\Psi^*_{\ep_1}(\b_1)\right)=
\ep_1 (2\h)^{1/2} \tilde G(\b)
  {\Ga{\fract{1}{2}+\fract{\b}{2\h}}\over
\Ga{\fract{\b}{2\h}}}
{\Ga{\fract{\b-\h}{\g}} \Ga{\fract{\g-\b-\h}{\g}} \over
\g\Ga{\fract{\g-2\h}{\g}}          }  \label{fin2II}
\eeq
We see that the vanishing  of this trace for $\g=2\h$ is due
to the function $\Ga{1-\fract{2\h}{\g}}$ in denominator of \r{fin2II}.

We can check that \r{fin2II} is in accordance with the normalization
condition for the type II intertwining operators 
\beq
\Psi^*_{\ep_2}(\b_2)\Psi^*_{\ep_1}(\b_1)={\ep_1 g \delta_{\ep_1,-\ep_2}
\over \b_1-\b_2-\h}+o(\b_1-\b_2-\h)
\label{normII}
\eeq
One can see that 
the function \r{fin2II} indeed has a pole at the point $\b_1=\b_2+\h$
with the residue equal to $\ep_1 g$.

\setcounter{equation}{0}

\section{Deformation of Gauss-Manin Connection}

As we have seen in the previous sections the general formula for the 
trace
simplifies much in the case when $\g=2\h$. In this case the
trace of the product of type II intertwining operators appear to 
be equals zero identically. 
To make a trace of type II vertex operators non-zero 
for $\g=2\h$ we should add some combination of type
I operators. In this and the next sections  
we are going to consider only simplest cases
of the formula \r{general-trace}, while general situation will
be considered elsewhere. Our goal is to identify some 
particular cases of the general trace formula
\r{general-trace} with those obtained by bootstrap approach
for $SU(2)$-invariant Thirring 
model. On the way of the identification we will find some identities 
between form factor integrals which turns out to be  deformations
of the Gauss-Manin connection identities for the hyperelliptic
curves. 

Let us start with the case of one type I
intertwining operator in the formula \r{general-trace} (which does not 
 correspond to local operators). 
For simplicity we consider this operator to be $\Phi_-(\z)$ component
of the intertwining operator. This means that $M=1$ and $m=0$ in
\r{general-trace}. The number $N$ and $n$ which describe the total 
number of type II operators and ``plus'' components of them respectively
are related because of the condition \r{condition}
\beq
N=2n+1\ .
\label{condI1}
\eeq
The case $N=1$ is trivial. The general formula \r{general-trace} does
not contain integration at all. First nontrivial case is $N=3$ and $n=1$.
In such a case the formula \r{general-trace} reads 
\bea 
&{1\over \tr\left(e^{2\h d}\right)}
\tr\left(e^{2\h d}
\prod_{j=1}^3\Psi^*_{\ep_j}(\b_j)\Phi_{-}(\z)
\right)
=2\h A_{3,1}^{1,0}\prod_{j<j'}\tilde\zeta(\b_j-\b_{j'})
\prod_{j}\overline\zeta(\z-\b_{j})\nn\\
&\qquad \times \int_{\tilde C}{d\v\over2\pi i}
\prod_{j=1}^{3}
\Ga{{\b_j-\v\over2\h}}\Ga{{\v-\b_j+\h\over2\h}}
{\tilde P_{\ep}(\v;\b)\over \prod_{j=1}^3 (\v-\b_j-\h)
}
\sin\pi\fract{\v-\z}{2\h}
\label{tr1,3}
\eea 
The contour $\tilde C$ in the integral \r{tr1,3} is shown on the Fig.~2.

\setlength{\unitlength}{0.00083333in}
\begin{picture}(6261,4756)(-1200,-310)
\thicklines
\put(100,   2100){\makebox(0,0)   [lb]{$\bullet$}}
\put(100,   1500){\makebox(0,0)   [lb]{$\bullet$}}
\put(100,    900){\makebox(0,0)   [lb]{$\bullet$}}
\put(1300,2100){\makebox(0,0)   [lb]{$\bullet$}}
\put(1300,1500){\makebox(0,0)   [lb]{$\bullet$}}
\put(1300, 900){\makebox(0,0)   [lb]{$\bullet$}}
\put(2500,2100){\makebox(0,0)   [lb]{$\bullet$}}
\put(2500,1500){\makebox(0,0)   [lb]{$\bullet$}}
\put(2500, 900){\makebox(0,0)   [lb]{$\bullet$}}
\put(3800,2100){\makebox(0,0)   [lb]{$\bullet$}}
\put(3800,1500){\makebox(0,0)   [lb]{$\bullet$}}
\put(3800, 900){\makebox(0,0)   [lb]{$\bullet$}}
\put(100,   2200){\makebox(0,0)   [cb]{$\b_1+2\h$}}
\put(100,   1600){\makebox(0,0)   [cb]{$\b_2+2\h$}}
\put(100,   1000){\makebox(0,0)   [cb]{$\b_3+2\h$}}
\put(1300,2200){\makebox(0,0)   [cb]{$\b_1+\h$}}
\put(1300,1600){\makebox(0,0)   [cb]{$\b_2+\h$}}
\put(1300,1000){\makebox(0,0)   [cb]{$\b_3+\h$}}
\put(2500,2200){\makebox(0,0)   [cb]{$\b_1$}}
\put(2500,1600){\makebox(0,0)   [cb]{$\b_2$}}
\put(2500,1000){\makebox(0,0)   [cb]{$\b_3$}}
\put(3800,2200){\makebox(0,0)   [cb]{$\b_1-\h$}}
\put(3800,1600){\makebox(0,0)   [cb]{$\b_2-\h$}}
\put(3800,1000){\makebox(0,0)   [cb]{$\b_3-\h$}}
\put(1300,100){\makebox(0,0)[cb]{$\tilde C$}}
\put(3800,100){\makebox(0,0)[cb]{$C'$}}
\put(2500,-200){\makebox(0,0) [cb]{Fig.~2.}}
\path(3375,12)(3375,4212)
\path(3405.000,4092.000)
     (3375.000,4212.000)
     (3345.000,4092.000)
\path(3075,8)	(3083.912,66.641)
	(3092.540,124.271)
	(3100.886,180.902)
	(3108.949,236.547)
	(3116.729,291.220)
	(3124.226,344.932)
	(3131.440,397.696)
	(3138.371,449.527)
	(3145.019,500.436)
	(3151.384,550.436)
	(3157.467,599.540)
	(3163.266,647.762)
	(3168.783,695.113)
	(3174.017,741.608)
	(3178.968,787.258)
	(3183.636,832.077)
	(3188.021,876.077)
	(3192.123,919.271)
	(3195.942,961.673)
	(3199.478,1003.295)
	(3202.732,1044.150)
	(3205.702,1084.251)
	(3208.390,1123.610)
	(3210.795,1162.241)
	(3212.916,1200.157)
	(3214.755,1237.370)
	(3217.584,1309.739)
	(3219.282,1379.452)
	(3219.847,1446.613)
	(3219.282,1511.323)
	(3217.584,1573.686)
	(3214.755,1633.804)
	(3210.795,1691.782)
	(3205.702,1747.722)
	(3199.478,1801.726)
	(3192.123,1853.899)
	(3183.636,1904.342)
	(3174.017,1953.159)
	(3163.266,2000.453)
	(3151.384,2046.327)
	(3138.371,2090.884)
	(3124.226,2134.227)
	(3108.949,2176.458)
	(3092.540,2217.682)
	(3075.000,2258.000)
\path(3075,2258)	(3056.519,2297.143)
	(3035.317,2338.700)
	(3011.435,2381.974)
	(2984.910,2426.265)
	(2955.784,2470.877)
	(2924.096,2515.110)
	(2889.885,2558.266)
	(2853.191,2599.648)
	(2814.054,2638.555)
	(2772.512,2674.292)
	(2728.607,2706.158)
	(2682.376,2733.456)
	(2633.861,2755.487)
	(2583.100,2771.554)
	(2530.133,2780.958)
	(2475.000,2783.000)
\path(2475,2783)	(2416.651,2776.834)
	(2361.217,2762.841)
	(2308.669,2741.802)
	(2258.978,2714.499)
	(2212.115,2681.715)
	(2168.052,2644.232)
	(2126.758,2602.832)
	(2088.206,2558.296)
	(2052.366,2511.408)
	(2019.210,2462.948)
	(1988.708,2413.700)
	(1960.831,2364.445)
	(1935.551,2315.966)
	(1912.839,2269.043)
	(1892.665,2224.461)
	(1875.000,2183.000)
\path(1875,2183)	(1854.919,2119.864)
	(1843.207,2050.198)
	(1840.107,2013.207)
	(1838.643,1974.930)
	(1838.661,1935.484)
	(1840.008,1894.985)
	(1842.533,1853.548)
	(1846.084,1811.290)
	(1850.506,1768.325)
	(1855.649,1724.770)
	(1861.360,1680.740)
	(1867.486,1636.352)
	(1873.874,1591.721)
	(1880.374,1546.963)
	(1886.831,1502.193)
	(1893.094,1457.527)
	(1899.010,1413.082)
	(1904.427,1368.972)
	(1909.192,1325.314)
	(1913.153,1282.224)
	(1916.157,1239.817)
	(1918.052,1198.209)
	(1918.686,1157.515)
	(1917.906,1117.852)
	(1915.560,1079.336)
	(1911.495,1042.081)
	(1897.600,971.822)
	(1875.000,908.000)
\path(1875,908)	(1856.413,868.967)
	(1834.922,827.709)
	(1810.609,784.882)
	(1783.556,741.138)
	(1753.845,697.133)
	(1721.557,653.521)
	(1686.775,610.957)
	(1649.580,570.095)
	(1610.055,531.589)
	(1568.280,496.094)
	(1524.340,464.264)
	(1478.314,436.754)
	(1430.285,414.218)
	(1380.335,397.311)
	(1328.546,386.687)
	(1275.000,383.000)
\path(1275,383)	(1221.454,386.687)
	(1169.665,397.311)
	(1119.715,414.218)
	(1071.686,436.754)
	(1025.660,464.264)
	(981.720,496.094)
	(939.945,531.589)
	(900.420,570.095)
	(863.225,610.957)
	(828.443,653.521)
	(796.155,697.133)
	(766.444,741.138)
	(739.391,784.882)
	(715.078,827.709)
	(693.587,868.967)
	(675.000,908.000)
\path(675,908)	(649.737,971.886)
	(628.565,1042.322)
	(619.464,1079.699)
	(611.328,1118.356)
	(604.137,1158.174)
	(597.873,1199.034)
	(592.515,1240.818)
	(588.044,1283.406)
	(584.441,1326.680)
	(581.687,1370.520)
	(579.761,1414.807)
	(578.646,1459.423)
	(578.321,1504.248)
	(578.768,1549.164)
	(579.966,1594.051)
	(581.896,1638.790)
	(584.540,1683.263)
	(587.878,1727.350)
	(591.889,1770.932)
	(596.557,1813.891)
	(601.859,1856.107)
	(607.778,1897.461)
	(614.295,1937.835)
	(621.389,1977.110)
	(637.232,2051.884)
	(655.154,2120.831)
	(675.000,2183.000)
\path(675,2183)	(704.206,2256.539)
	(721.699,2295.111)
	(741.000,2334.726)
	(762.018,2375.264)
	(784.664,2416.606)
	(808.846,2458.633)
	(834.476,2501.225)
	(861.464,2544.264)
	(889.719,2587.629)
	(919.152,2631.202)
	(949.672,2674.864)
	(981.191,2718.495)
	(1013.617,2761.976)
	(1046.861,2805.187)
	(1080.832,2848.010)
	(1115.442,2890.325)
	(1150.601,2932.013)
	(1186.217,2972.954)
	(1222.201,3013.030)
	(1258.464,3052.120)
	(1294.915,3090.107)
	(1368.023,3162.290)
	(1440.804,3228.625)
	(1512.541,3288.158)
	(1582.512,3339.934)
	(1650.000,3383.000)
\path(1650,3383)	(1713.030,3411.285)
	(1784.476,3429.363)
	(1822.982,3435.131)
	(1863.144,3439.015)
	(1904.813,3441.238)
	(1947.840,3442.022)
	(1992.075,3441.589)
	(2037.370,3440.164)
	(2083.576,3437.967)
	(2130.543,3435.222)
	(2178.122,3432.151)
	(2226.163,3428.978)
	(2274.519,3425.923)
	(2323.039,3423.211)
	(2371.574,3421.064)
	(2419.975,3419.704)
	(2468.094,3419.354)
	(2515.780,3420.237)
	(2562.885,3422.574)
	(2609.259,3426.590)
	(2654.753,3432.506)
	(2699.219,3440.545)
	(2742.506,3450.930)
	(2784.466,3463.883)
	(2824.950,3479.627)
	(2863.808,3498.384)
	(2936.050,3545.829)
	(3000.000,3608.000)
\path(3000,3608)	(3029.909,3649.715)
	(3051.272,3695.185)
	(3064.090,3745.949)
	(3068.362,3803.544)
	(3064.090,3869.508)
	(3058.749,3906.110)
	(3051.272,3945.380)
	(3041.658,3987.513)
	(3029.909,4032.698)
	(3016.022,4081.130)
	(3000.000,4133.000)
\path(3065.006,4027.766)(3000.000,4133.000)(3007.837,4009.555)
\end{picture}

The  form of the contour $\tilde C$ and presence of
``strange'' poles (from bootstrap approach point of view) 
at the points $\b_j+\h$ in the integral \r{tr1,3} were  main obstacles
in identification of free filed \cite{JM} and bootstrap \cite{S1}
formulas for form factors. We would like to state that both
these obstacles can be overcome because of the following 
\medskip

\noindent {\bf Proposition 3.} {\it For} $m=1,2,3$  
{\it we have identities}  
\bea
&\prod_{j\neq m}(\b_m+\h-\b_j)
\int_{\tilde C}{d\v\over 2\pi i}
\prod_{j=1}^3\varphi(\v-\b_j)
{1\over \v-\b_m-\h}\exp\left({\pm i\frac{\pi\v}{2\h}}\right)\nn\\
&\qquad =
\int_{C'}{d\v\over 2\pi i}
\left[\prod_{j=1}^3\varphi(\v-\b_j)\right]
 (\v-\b_m)\exp\left({\pm i\fract{\pi\v}{2\h}}\right)
\label{def-ident}
\eea
{\it where}
\beq
\varphi(\v)=\Ga{-\frac{\v}{2\h}}\Ga{\frac{1}{2}+\frac{\v}{2\h}}
\label{varphi}
\eeq
{\it which satisfy}
\beq
\varphi(\v+2\h)=-{\v+\h\over\v+2\h}\varphi(\v)
\label{quasi}
\eeq
{\it and the contour $C'$ is a straight line between the points $\b_j$ and
$\b_j-\h$, $j=1,2,3$ and shown on the Fig.~2.}
\smallskip

The identity formulated in \r{def-ident} is of a type  
of the total difference identities found by F.A.~Smirnov in 
\cite{S3,S4}. The proof of \r{def-ident} is based on the
property \r{quasi} of the function $\varphi(\v)$.

Hence due to \r{def-ident} 
the extremal  component ($\ep_3,\ep_2,\ep_1={+}{-}{-}$) 
of the three-particle form factor of the operator $\Phi_-(\z)$
is equal to
\bea 
&{1\over \tr\left(e^{2\h d}\right)}
\tr\left(e^{2\h d}
\Psi^*_{+}(\b_3)\Psi^*_{-}(\b_2)\Psi^*_{-}(\b_1)
\Phi_{-}(\z)
\right)\nn\\
&\quad
=4\h A_{3,1}^{1,0}\prod_{j<j'}\tilde\zeta(\b_j-\b_{j'})
\prod_{j}\overline\zeta(\z-\b_{j})
\prod_{j=1}^2{1\over(\b_3+\h-\b_j)}
\nn\\
&\qquad \times
\exp\left({i\pi\over4\h}\left(3\z+
\sum_{j=1}^{3}\b_j\right)\right)
\left[
\prod_{j=1}^{3}\sin\pi\fract{\z-\b_j}{4\h}+i
\prod_{j=1}^{3}\cos\pi\fract{\z-\b_j}{4\h}
\right]
\nn\\
&\qquad \times 
\int_{C'}{d\v\over 2\pi i}
\left[\prod_{j=1}^3\varphi(\v-\b_j)\right]
 (\v-\b_3)\exp\left({ -\fract{i\pi\v}{2\h}}\right)
\label{tr1,3fin}
\eea 
The formula \r{tr1,3fin} coincides with
corresponding formula
from \cite{JKMQ}. 
To obtain \r{tr1,3fin} an obvious identity
\beq
\int_{C'}{d\v\over 2\pi i}
\left[\prod_{j=1}^3\varphi(\v-\b_j)\right]
 (\v-\b_3)
\left[
\exp\left({ \fract{i\pi\v}{2\h}}\right)
+B_2\exp\left({ -\fract{i\pi\v}{2\h}}\right)
\right]=0
\label{obvious}
\eeq
have been used, where
$$
B_2=
\exp\left(\sum_{j=1}^{3} \frac{i\pi\b_j}{2\h}\right)
\sum_{j=1}^{3}
\exp\left(- \frac{i\pi\b_j}{2\h}\right)\ .
$$

Another important application of the identities \r{def-ident}
is that they provide a deformation of Gauss-Manin connection
for elliptic curves \cite{GM}. To explain what we mean let us define
two integrals
\beq
F^{(k)}_p(\b_1,\b_2,\b_3)
=\int_{C'}{v^p\, d\v\over2\pi i}
\prod_{j=1}^3\varphi(\v-\b_j) \exp\left({ik\fract{\pi\v}{2\h}}
\right),\quad
p=0,1\quad k=\pm 
\label{def-hyp}
\eeq
which are deformed analogs of first and second kind integrals
on elliptic curve \cite{S3}. Due to \r{def-ident} the integrals
\r{def-hyp} 
satisfy the difference equation 
\beq
F^{(k)}_p(\b_1,\b_2,\b_3+2\h)+\sum_{r=1}^2
A_{pr}(\b_1,\b_2,\b_3)
F^{(k)}_r(\b_1,\b_2,\b_3)=0
\label{dif-eq}
\eeq
with the connection matrix 
\beq
A(\b)=
\left(
\begin{array}{cc}
1&0\\ \h&1
\end{array}
\right)
+
{1\over\prod_{j=1}^2(\b_3+\h-\b_j)}
\left(
\begin{array}{cc}
\h\b_3&\h\\ \h\b_3(\b_3+\h)&\h(\b_3+\h)
\end{array}
\right)
\label{def-con}
\eeq
Difference equation
\r{dif-eq} is a deformed analog of the Gauss-Manin connection.
Classical limit ($\h\to 0$) of the integrals 
$F^{(k)}_p(\b_1,\b_2,\b_3)$ is not straightforward and was
investigated in \cite{S3,S4}. In such limit these integrals
go to elliptic integrals
\bea
F^{(+)}_p(\b_1,\b_2,\b_3)&
\stackreb{\sim}{\h\to 0}& 
\int_{\b_1}^{\b_2}
{v^p\ d\v\over
\sqrt{\prod_{j=1}^3(\v-\b_j)}}\nn\\
F^{(-)}_p(\b_1,\b_2,\b_3)&
\stackreb{\sim}{\h\to 0}& 
\int_{\b_2}^{\b_3}
{v^p\ d\v\over
\sqrt{\prod_{j=1}^3(\v-\b_j)}}
\nn
\eea 
and difference equation \r{dif-eq} to the
classical Gauss-Manin connection \cite{GM}.
Identity \r{def-ident} becomes the identity between elliptic integrals 
of the type
\beq
\prod_{j\neq m}(\b_3-\b_j)
\int_{\b_1}^{\b_2} 
{d\v\over 
\sqrt{\prod_{j=1}^3(\v-\b_j)}
(\v-\b_3)}=
\int_{\b_1}^{\b_2} 
{(\v-\b_3)d\v\over 
\sqrt{\prod_{j=1}^3(\v-\b_j)}}
\label{cl-ident}
\eeq
which describe a relation between second kind elliptic integrals
with singularities at the point $\b_3$ and at the point $\infty$.



The next example which is interesting from physical point of view 
is two type I intertwining operators
in the trace formula \r{general-trace}.
The case when $m=n=0$ and $M=N=2$ 
is trivial because in this
case \r{general-trace} does not contain integrations.
We can state that  the integral \r{general-trace} can be reduced
to form factor integrals in $SU(2)$-invariant Thirring model for
$N=2n$, $M=2m=2$, $\g=2\h$
 and arguments of type I operators 
are related  as
\beq
\z_1=\z_2+\h\ .
\label{phys}
\eeq 
Because of the property \r{inverse} 
this case correspond to $2n$-particle
form factor of the generating function of the local operators \cite{L}
\beq
\Lambda(\z)=\Phi_{+}(\z)
\Phi_{-}(\z+\h) 
+\Phi_{-}(\z) \Phi_{+}(\z+\h) =
\Phi_{-}(\z) \Phi_{-}(\z+\h) f_0-f_0\Phi_{-}(\z) \Phi_{-}(\z+\h) 
\label{loc-Luk}
\eeq 
The integral formula \r{general-trace} reduces in this case to the
multiple integral 
\bea 
&{1\over\tr\left(e^{2\h d}\right)}\tr\left(e^{2\h d}
\prod_{j=1}^{2n}\Psi^*_{\ep_j}(\b_j)\Phi_{+}(\z)\Phi_{-}(\z+\h)
\right)\nn\\
&\quad =  {A_{2n,2}^{n,1} \pi^{3/2}(-1)^n 
\over
 2^{n+1}  (2\h)^{n(n-3)/2}  
G^{1/2}(0) }
{\dis \prod_{j<j'}\tilde\zeta(\b_j-\b_{j'})\over\dis
\prod_{j=1}^{2n}\cos\pi\fract{\z-\b_{j}}{2\h}}\nn\\
&\qquad  \prod_{k=1}^n \int_{\tilde C}{d\v_k\over2\pi i}
\prod_{j=1}^{2n}
{\varphi(\v_k-\b_j)\over
(\v_k-\b_j-\h)}
\tilde P_{\ep}(\v;\b)
\prod_{k<k'} (v_k-v_{k'}-\h)\nn\\
&\qquad\times
\prod_{k<k'}
\sin\pi\frac{\v_k-\v_{k'}}{\h}
\prod_{k=1}^n 
\sin\pi\frac{\v_k-\z}{\h} \int_{C}{d\u\over2\pi i}
{\dis
\prod_{j=1}^{2n}\sin\pi\frac{\u-\b_j}{2\h}
\over
\dis
\prod_{k=1}^n\sin\pi\frac{\u-\v_k}{\h}
\sin\pi\frac{\u-\z}{\h}     }
\label{mix2-2n}
\eea

Note that in order to set $\z_1=\z_2+\h$ 
in the general formula \r{general-trace}
we have used the explicit form of the polynomial 
$P_{+-}(u)$ which cancel the pole
of the function $\Ga{{\u-\z_1\over2\h}}$ at the point $\z_1$. It means that
there is no pinching of the contour $C$ between points $\z_2+\h$ and
$\z_1$ when $\z_1\to\z_2+\h$.
The poles of the integrand over $u$ at the points 
$$\z+r\h,\quad \v+r\h,\quad k=1,\ldots,n $$ 
are inside the contour $C$ for $r=1,2,\ldots$ 
and are outside  $C$ for  $r=0,{-}1,{-}2,\ldots\,$.
Note also that all $\Gamma$-function depending on the integration
variable $u$ turns into sin-functions, so we can calculate the 
integral over $u$.  
After this  
we can see
that one more integration
over variables $\v_k$ can be calculated and 
the result  of calculation  
coincides with
form factors integral obtained in \cite{KS,S1} for $SU(2)$-invariant
Thirring model.

Unfortunately, this treatment of the  integral \r{mix2-2n}
is technically complicated 
problem. We are going to publish
this calculation  as a separate paper \cite{KLP2}. 
Here we would like to point out the main steps of our 
construction.

First step is the calculation of the integral over $u$ in \r{mix2-2n}
which can be done easily. Second one is the calculation
of the integral over one of variables $v_k$ using 
Smirnov's hint (see section 6). 
After these calculations we left with $n-1$
integrals over the contours $\tilde C$ of the type shown
on Fig.~2. In order to rectify contours
we should use the identities
which are analogs of  
the Gauss-Manin connection identities.
In the simplest case $n=2$ this identity
reads
\bea
&\int_{\tilde C}{d\v\over 2\pi i}
{\dis \prod_{j=1}^4\varphi(\v-\b_j)\over \v-\b_4-\h}
\exp\pi\left(\pm {i\pi\v\over\h}  \right)                         \nn\\
&\qquad ={1\over\prod_{j\neq 4}(\b_4-\b_j+\h)}
\int_{C'}{d\v\over 2\pi i}
d(\v,\b)\prod_{j=1}^4\varphi(\v-\b_j)
\exp\left(\pm{i\pi\v\over\h}\right)
\label{Gauss-Manin}\\
&d(\v,\b)=2\v^2-\v\left(2\h-\sum_{j=1}^4\b_j\right)+
\b_4(\b_1+\b_2+\b_3-\b_4-2\h)\nn
\eea
The contours $\tilde C$ and  $C'$ in \r{Gauss-Manin} 
are the same as shown on the Fig.~2.
These identities 
can be easily generalized to arbitrary $n\geq3$ and 
are the  analogs of 
the Gauss-Manin connection identities for the hyperelliptic
curves of the genus $n-1$. 
The last step of identification \r{mix2-2n} with form factor
 integrals is using the identities of the type \r{obvious}.
The identities of this kind  allow one 
to extract from  \r{mix2-2n} the integrals which do not  
depend on the
spectral parameter $\z$ of the generating function $\Lambda(\z)$.
Realization of such a program 
for the trace 
$$ 
{\dis\tr\left(e^{2\h d}
\Psi^*_{+}(\b_4)\Psi^*_{+}(\b_3)
\Psi^*_{-}(\b_2)\Psi^*_{-}(\b_1)
\Lambda(\z)
\right)\over\dis
\tr\left(e^{2\h d}
\right)}
$$
leads to the integral
\bea
&\prod_{j<j'}\tilde\zeta(\b_j-\b_{j'})
\prod_{i=1,2\atop j=3,4}{1\over (\b_j-\b_i+\h)}\nn\\
&\quad \times
\left(\sum_{j=1}^4\exp{i\pi\b_j\over\h}\right)^{-1}
\exp\left(\sum_{j=1}^4{i\pi\b_j\over2\h}\right)\nn\\
&\quad \times\int_{C'}{d\v\over2\pi i}
\prod_{j=1}^{4} \varphi(\b_j-\v)\exp\left(-{i\pi v\over\h}\right)
\left[\prod_{j=1,2}(\v-\b_j+\h)+\prod_{J=3,4}(\v-\b_j)
\right]
\label{Smir-ff}
\eea
Formula
\r{Smir-ff} coincide with analogous formula from \cite{S2}.

\newpage
  
\section{Acknowledgements}

Authors would like to thank  
V.~Fock, A.~Gerasimov, A.~Levin and M.~Olshanetskii
for helpfull discussions.
First author (S.K.) was supported by RFBR grant 96-01-01106 and 
INTAS grant 1010-CT93-0023. The second author (D.L.) was supported 
by RFBR grant 96-01-01101, INTAS grant 1010-CT93-0023, 
Wolkswagen Stiffung and AMSFSU grant. 
D.L. would like to acknowledge O.~Lechtenfeld
and S.~Ketov for warm hospitality at
Institut f\"ur Theoretische Physik, Universit\"at Hannover,
where part of this work was done.
Last author (S.P.) 
would like to acknowledge 
Physikalisches Institut, Universit\"at Bonn,
Max-Planck Institute f\"ur Mathematik (Bonn) and 
Institut f\"ur Theoretische Physik, Universit\"at Hannover
for the hospitality.
The research of the last author (S.P.) 
was supported by Heisenberg-Landau 
program and  RFBR grant 96-01-01106.

\end{document}